\documentclass[11pt,twoside]{article}
\usepackage{graphicx}
\usepackage{rotating}
\usepackage{amsmath}
\usepackage{amsfonts}
\usepackage{amssymb}
\usepackage{amsthm}
\usepackage{tabularx}
\usepackage[all]{xy}

\newtheorem{thm}{Theorem}[section]
\newtheorem{pro}[thm]{Proposition}
 \newtheorem{cor}[thm]{Corollary}
 
 \newtheorem{rem}[thm]{Remark}
 \newtheorem{lem}[thm]{Lemma}
  \newtheorem{defn}[thm]{Definition}

\begin{document}


 \centerline{\large{\textbf{   ON CLASSIFYING   HUREWICZ  FIBRATIONS  }}}
\centerline{\large{\textbf{    AND FIBRE  BUNDLES OVER  }}}
\centerline{\large{\textbf{     POLYHEDRON BASES }}}
 \centerline{\large{\textbf{     }}}
\centerline{}

\centerline{{\emph{\textbf{Amin  Saif and Adem K{\i}l{\i}\c cman}}}}
\centerline{\footnotesize{Institute of Mathematical Research
(INSPEM), University Putra Malaysia}}
\centerline{\footnotesize{43400 UPM, Serdang, Selangor,
Malaysia}}


\begin{abstract}
\footnotesize{ Let $f:E\longrightarrow O$ be a Hurewicz fibration
with a fiber space $F_{r_{o}}$ and a lifting function $L_{f}$. The
\emph{$Lf-$function} $\Theta_{L_{f}}$  of $f$ is defined by the
restriction map of $L_{f}$ on the  space $\Omega(O,r_{o})\times
F_{r_{o}}\times \{1\}$.  The purpose of this paper is to give some
results which show  the role of $Lf-$functions in finding a fiber
homotopically equivalent relation between two fibrations, over a
common polyhedron base. Furthermore we will prove the equivalently
between our results and Dold's theorem
in fiber bundles, over a common suspension base of polyhedron spaces.   
}

\end{abstract}


\section{Introduction}
In what follows, for a topological space $E$, $\widetilde{E}$ will
denote a set of all constant path $\widetilde{e}$ into $e\in E$,
$\overline{\alpha}$ the inverse path of $\alpha\in E^I$, $\star$
the usual path multiplication operation, $ \simeq  $ the same
homotopy type for spaces and homotopic for maps and $\simeq_{f}$ a
fiber homotopy.

First we recall from \cite{Fedorchuk} that the  \emph{simplicial
complex} $K$ contains of a set $\{v\}$ of vertices and a set
$\{s\}$ of finite nonempty subsets of $\{v\}$ is called simplexes
such that any set consisting of exactly one vertex in a  simplex
and any nonempty subset of simplex is a simplex.  $/K/$ denotes a
set of all functions $\beta$ from the  set of vertices of $K$ to
$I$ such that for any function $\beta  \in /K/$, the set $\{v\in
K: \beta(v)=0\}$ is a simplex and  $\sum_{v\in K} \beta(v)=1$.
  The topology  on $/K/$  is a topology which
  induced by metric $d$ on $/K/$ defined by
  $d(\alpha,\beta)=\{\sum_{v\in K} [\alpha(v)-\beta(v)]^2\}^{0.5}$.
A topological space $E$ is called a \emph{polyhedron} if there is
a simplicial complex $K$ and a homeomorphism $f:/K/\longrightarrow
E$. A closed subspace $A$ of $E$ is called a \emph{subpolyhedron}
of $E$ if there is simplicial complex $L\subset K$ such that
$f(/L/)=A$.
 And if $A$ is a \emph{subpolyhedron} of $E$, we say that the pair $(E,A)$
  \emph{polyhedron pair}. In our paper, for any simplicial
complex, any vertex belong to a finite simplexes, hence any
polyhedron in this case will be an ANR (see \cite{Fedorchuk}).

From \cite{Alexander} the \emph{suspension } $S(O)$  of $O$ based
a fixed point $r_o\in O$ is defined to be the quotient space of
$O\times I$ in which for all $b\in O$, $(b,0)$ is identified to
$(r_o,0)$ and $(b,1)$ is identified to $(r_o,1)$. Consider $S(O)$
as the union of two cones, one of them is defined
\[S_{0}(O)=\{[(x,t)]\in S(O):  (x,t)\in O\times [0,1/2]\}\] with
$(x,0)$ identified to $(r_o,0)$ and
       the other  by \[S_{1}(O)=\{[(x,t)]\in S(O):  (x,t)\in O\times [1/2,1]\}\] with $(x,1)$ identified to $(r_o,1)$.
        The cones are always contractible spaces (see \cite{KarimovD1, KarimovD2}).
       If $ O$ is an absolute neighborhood retract (ANR)  space, then $S(O)$ is  an ANR.  Similarly, for
       polyhedron property, for more details see \cite{Klein, Stramaccia}.

Dold's theorem is one of the famous solutions for the
classification problem in fiber bundles over  the    $n-$sphere
bases $S^{n}$, (see\cite{Husemoller}). James in \cite{James}
showed that Dold's theorem remains valid if we use suspensions of
polyhedron space  instead of $n-$spheres $S^n$ for the base of
bundles., that is, Dold's theorem will take the following form:

\begin{thm}\label{Am1}
\emph{[The Dold's theorem]  Let $\gamma =(E ,f ,S(O),F ,G )$ and
$\gamma'=(E',f',S^{n},F',G')$ be two fiber bundles over a
suspension $S(O)$ of polyhedron space $O$ with locally compact
fibers $F $and $F'$. Let $\mu :(O, x_{o})\longrightarrow (G ,e )$
and $\mu':(O, x_{o})\longrightarrow (G',e')$ be characteristic
maps of $\gamma $ and $\gamma'$, respectively and let $i :G
\longrightarrow H(F,F)$ and $i':G'\longrightarrow H(F', F')$ be
the inclusion maps. Then $\gamma $ and $\gamma'$ are fiber
homotopy equivalent if and only if there is homotopy equivalence
$g: F \longrightarrow F'$ such that the maps \[ q (x)=  g\circ (i
\circ \mu )(x)\circ \overleftarrow{g}\quad \mbox{ and}\quad
q'(x)=(i'\circ \mu')(x)\] from $O$ into $H(F',F')$ are homotopic,
where $H(F,F)$ is the set of all homotopy equivalences from $F$
into $F$ and $\overleftarrow{g}$ is the inverse homotopy of $g\in
H(F,F)$.}
\end{thm}

 This paper is organized as follows. It consists of four
sections. After this Introduction, Section 2 is devoted to some
preliminaries. In Section 3 we shall start  by giving some results
about homotopy extension property and $Lf-$function properties.
Next we show the role of $Lf-$ function $\Theta_{L_{f}}$ in
finding fiber homotopy equivalence between two fibrations. Mainly,
we prove the following theorem:
\begin{thm}\label{Am2}
Let $[E_1,f_1,O,F^1_{r_{o}}]$ and $[E_2,f_2,O,F^2_{r_{o}}]$ be two
fibrations over a polyhedron base $O$. Let $O$ be  the union of
two subpolyhedra   $ O_{1} $ and
  $ O_{2} $ such that $O_{1}$ is  a contractible in
$O$   to   $r_{o}\in O_{3}=O_{1}\cap O_{2}$ leaves $r_{o}$ fixed
and $ O_{2} $ is    contractible to $r_{o}$. If $ O_{3}$ is
subpolyhedra  of $O$, then $f_{1}$ and $f_{2}$ are fiber homotopy
equivalent if and only if they have conjugate $Lf-$functions by
$g\in H(F^{1}_{r_{o}},F^{2}_{r_{o}})$.
\end{thm}
In section 4 we will apply $Lf-$function in fiber bundles by
proving  the equivalently between Theorem \ref{Am2} and Dold's
theorem.

All topological spaces in  this paper will be assumed Hausdorff
spaces.

\section{Preliminaries }

Recall \cite{Allen} that the path space $Pa(E,e_{o})= \{\alpha\in
E^I: \alpha(0)=e_{o}\}$
 based at fixed point $e_{o}$, a loop space $\Omega(E,e_{o})$ and $\widetilde{E}$
 are closed in a path space $E^I$. If $E$ is a metrizable (\emph{resp.} ANR),
 then $E^I$, $Pa(E,e_{o})$ and $\Omega(E,e_{o})$ are  metrizable (\emph{resp.}
 ANR).

 Let $f:E\longrightarrow O$ be a fibration with a base $O$, total space
$E$ and fiber space $F_{r_{o}}=f^{-1}(r_{o})$, where $r_{o}\in O$.
A map $L_{f}: \bigtriangleup f \longrightarrow E^I$ is called a
\emph{ lifting function} for $f$ if $L_{f}(e,\alpha)(0)=e$ and
$f[L_{f}(e,\alpha)]=\alpha$ for all $(e,\alpha)\in \bigtriangleup
f$, where  $\bigtriangleup f = \{(e,\alpha)\in E\times
Pa(O):f(e)=\alpha(0)\}$. If $L_{f}(e,f\circ
\widetilde{s})=\widetilde{e} \mbox { } \mbox{ for all } \mbox{
}e\in E$, then   the  lifting function is called a \emph{  regular
lifting function}. A  fibration $f$ is called \emph{regular
fibration} if it has regular lifting function (see \cite{Allen}).

Curtis-Hurewicz theorem, \cite{Hur1}, is one of the famous
theorems  in  fibration theory which shows that any map is regular
fibration if and only if it has regular lifting function.


\begin{thm}\label{Am3}
\emph{[Dold-Fadell Theorem] (see \cite{Fad1}) Let
$f_{1}:E_{1}\longrightarrow O$ and $f_{2}:E_{2}\longrightarrow O$
be two fibrations over  an ANR pathwise connected base $O$. Then
 $f_{1}$ and $f_{2}$ are
fiber homotopy equivalent if and only if there is a fiber map
$h:E_{1}\longrightarrow E_{2}$ such that the restriction map of
$h$ on $f^{-1}_{1}(r_{o})$ is homotopy equivalence into $
f^{-1}_{2}(r_{o})$, for some $r_{o}\in O$.}
\end{thm}

  A closed subspace $A$ of a space $E$ is said to have a
\emph{homotopy extension property} in $E$ with respect to a space
$O$ if any map $f:(E\times \{0\})\cup (A\times I)\longrightarrow
O$ can be extended to a map $F:E\times I\longrightarrow O$.

\begin{thm}\label{Am4}
\emph{(See\cite{Schirmer}  For a polyhedron pair $(E,A)$, $A$ has
a homotopy extension property in $E$ with respect to any space
$O$. }
\end{thm}

\begin{thm}\label{Am5}
\emph{(See \cite{Stramaccia}) Let $A$ be a closed subspace of
metrizable space $E$ and $O$ be an ANR space. Then  $A$ has a
homotopy extension property in $E$ with respect to a space $O$. }
\end{thm}
\begin{thm}\label{Am6}
\emph{ (See \cite{Stramaccia}) An ANR  closed subspace $A$ of an
ANR space $E$  has a homotopy extension property in $E$ with
respect to any space $O$. }
\end{thm}

Recall \cite{JamesA} that if  $f:E\longrightarrow O$ is  a
fibration and $A$ is    a subspace of $O$ then the restriction map
$f|_{f^{-1}(A)}:f^{-1}(A) \longrightarrow A$ of $f$ is a fibration
and we denote  it by $f|A$.


\begin{thm}\label{Am7}
 \emph{(see \cite{JamesA}) Let $f_{1}:E_{1}\longrightarrow O$ and
$f_{2}:E_{2}\longrightarrow O$ be two
 regular fibrations over a polyhedron base $O$ and  $A$ be a subpolyhedron   of $O$. Suppose that  there are two fiber maps
$h_{1},h_{2}:f_{1}|A\longrightarrow f_{2}|A$  such that
$h_{1}\simeq_{f} h_{2}$. Then if  $k_{1}$ has an extension fiber
map $H_{1}:E_{1}\longrightarrow E_{2}$,   $h_{2}$ has an extension
 fiber map $H_{2}:E_{1}\longrightarrow E_{2}$  and $H_{1}\simeq_{f} H_{2}$.}
\end{thm}

Here we recall the details of the definition of fiber bundle which
will be used in the our results in Section 4.

\begin{defn}\label{Am8}
\emph{\cite{Husemoller} Let $E$, $O$  and $F$ be spaces. Let
$f:E\longrightarrow O$ be a map of $E$ onto $O$ and $G$ be group
of all homeomorphisms of $F$ onto $F$ with as a binary usual
composition operation $\circ$. Then $\gamma=(E,f,O,F,G)$ is said
to be \emph{ a fiber bundle over a base $O$} if there is an open
covering  $\{V_j:j\in \wedge\}$ of $O$ (where $\wedge$ is an index
set) and for each $j\in \wedge$, there is a homeomorphism $
\theta_j:V_j \times F\longrightarrow f^{-1}(V_j)$ such that:
\begin{enumerate}
    \item  $f[\theta_j(b,y)]=b$ for all $b\in V_j, y\in
    F$.
    \item For each pair $i,j\in \wedge$ and $b\in V_i\cap V_j$, the homeomorphism
    $\theta^{-1}_{jb}\circ \theta_{ib}:F\longrightarrow F$ corresponds to an element of $G$, where
     $\theta_{kb}: F\longrightarrow f^{-1}(b)$ defined by
     $\theta_{kb}(y)=\theta_k(b,y)$ for all $b\in V_k, y\in
    F$, ($k=i,j$).
    \item For each pair $i,j\in \wedge$, the function
    $g_{ij}:V_i\cap V_j\longrightarrow G$ given by $g_{ij}(b)=\theta^{-1}_{jb}\circ
    \theta_{ib}$ is a map.
\end{enumerate}}
\end{defn}
\begin{rem}\label{Am9}
 \emph{In fiber bundle $\gamma=(E,f,O,F,G)$, the maps
$\theta_j:V_j \times F\longrightarrow f^{-1}(V_j)$ are called the
\emph{coordinate functions}, the maps
$g_{ij}(b)=\theta^{-1}_{jb}\circ
    \theta_{ib}$ are called \emph{coordinate transformations}, the space $E$ is called \emph{bundle }over base $O$, and $F$ is called
   \emph{ fiber of bundle} $E$. We shall denote the identity element of a group $G$ by $\mathbf{g}$,
        the inverse element $g\in G$ by  $g^{-1}$.}
\end{rem}

\begin{thm}\label{Am10}
 \emph{\cite{James} Let $S^n$ be the $n-$sphere in $R^{n+1}$.
 For a fiber bundle $\gamma=(E,f,S^n,F,G)$, there is a characteristic map  $\mu:(S^{n-1}, x_o)\longrightarrow (G,e)$,
  , where $n>0$ is a positive integer.}
      \end{thm}

\section{An $Lf-$functions of fibration}
In this section, we will define the $Lf-$function and study its
properties. Next we will show its role in finding a fiber
homotopically equivalent relation between two fibrations over a
common polyhedron base.

\begin{thm}\label{Am11}
   Let $f:E\longrightarrow O$ be a
 regular fibration. Let  $(X,A)$ be an ANR pair (resp. be
 Polyhedron pair). If there is a map
   $G:(X\times \{0\})\cup (A\times I)\longrightarrow E$
 such
that
\[f[G(a,t)]=f[G(a,0)] \quad \mbox{ for  } \mbox{ } a\in A, t\in I,\]
then there is    a map $H:X\times I\longrightarrow E$ such that
$H$ is an extension of $G$,\[f[H(x,t)]=f[H(x,0)] \quad \mbox{ for
 } \mbox{ } x\in X, t\in I.\]
\end{thm}
\noindent \textbf{Proof.} By Theorem \ref{Am6} if  $(X,A)$ is  an
ANR pair or by Theorem \ref{Am4} if $(X,A)$ is a Polyhedron pair
we get that  the map $G$ can be extended to
 a map
$F :X\times I\longrightarrow E$. For a path  $\alpha\in E^I$ and
$r\in I$, we can define the path $\alpha_{r}$ in $E^I$ by
$\alpha_{r}(t)=\alpha[(1-t)r] $ for all $ t\in I$.  Hence we can
define the map $H:X\times I\longrightarrow E$ by
\[H(x,t)=L_{f}[F(x,t),f\circ F(x)_{t}](1)\quad \mbox{ for  } \mbox{ } x\in X, t\in I.\]
 At case $X\times \{0\}$,  by the
 regularity of $L_{f}$ and since $F$ is an extension
for $G$, we observe that for $x\in X$,
\begin{eqnarray*}
  H(x,0) &=& L_{f}(F(x,0),f\circ F(x)_{0})(1) \\
      &=&L_{f}(F(x,0),f\circ \widetilde{F(x,0)})(1) \\
   &=& [\widetilde{F(x,0)}](1)= F(x,0)= G(x,0).
\end{eqnarray*}
At case $A\times I$, since for $a\in A$ and $r,t\in I$,
\begin{eqnarray*}
  (f\circ F(a)_{t})(r) &=& f[ F(a)((1-r)t)]= f[F(a,(1-r)t)] \\
   &=&f[ G(a,(1-r)t)]= f[G(a,0)]= f[G(a,t)]\\
   &=&f[F(a,t)]=[f\circ \widetilde{F(ax,t)}](r),
   \end{eqnarray*}
   then  by the   regularity   of $L_{f}$ we get that $H(a,t)=F(a,t)=G(a,t)$ for all $t\in I,a\in A$. That is, $H$ is
   an extension for $G$.

 Finally, we also observe that
\begin{eqnarray*}
  f[H(x,t)]&=& f[L_{f}(F(x,t),f\circ F(x)_{t})(1)]  =(f\circ F(x)_{t})(1) \\
     &=& f[ F(x,0)]= f[ G(x,0)]=f[H(x,0)],
\end{eqnarray*}
for all  $x\in X$,  $t\in I$.\quad $\square$

We can give another rephrasing of  Theorem above in the following
corollary:

\begin{cor}\label{Am12}
 Let $f_{1}:E_{1}\longrightarrow O$ and
$f_{2}:E_{2}\longrightarrow O$ be two
 regular fibrations. Let $A$ be a closed subspace of $O$ and   $(E_{1},f_{1}^{-1}(A))$ be an ANR pair (resp. be
 Polyhedron pair). If there are two fiber maps
$k_{1},k_{2}:f_{1}|A\longrightarrow f_{2}|A$  such that
$k_{1}\simeq_{f} k_{2}$ and   $k_{1}$ has an extension fiber map
$K_{1}:E_{1}\longrightarrow E_{2}$, then  $k_{2}$ has an extension
 fiber map $K_{2}:E_{1}\longrightarrow E_{2}$  and $K_{1}\simeq_{f} K_{2}$.
\end{cor}
\noindent \textbf{Proof.}
  Since $k_{1}\simeq_{f} k_{2}$ then there is a homotopy $R:
  f_{1}^{-1}(A)\times I\longrightarrow f_{2}^{-1}(A)$ between
  $k_{1}:=R_{0}$ and $k_{2}:=R_{1}$ such that
  $f_{2}[R(a,t)]=f_{1}(a)$.
  We can apply Theorem \ref{Am11} on the regular fibration $f_{2}$ by taking
  $(X,A):=(E_{1},f_{1}^{-1}(A))$ and
  \[G(e,t)=\left\{
  \begin{array}{c l}
    R(e,t) & \mbox{for} \quad (e,t)\in f_{1}^{-1}(A)\times I\\
   K_{1}(e)   & \mbox{for} \quad (e,t)\in E_{1}\times \{0\},
  \end{array}
\right. \] to get the extension homotopy $H$ of $G$. Hence
$H_{0}=K_{1}$ and take $K_{2}=H_{1}$. Then $H$ is a homotopy
between two fiber maps $K_{1}$ and $K_{2}$ and we get \[
f_{2}[H(e,t)]  =f_{2}[H(e,0)]=f_{2}[K_{1}(e)]=f_{1}(e)\]for all
$e\in E_{1}$,  $t\in I$. That is, $K_{1}\simeq_{f} K_{2}$\quad
$\square$

\begin{rem}\label{Am13}
 \emph{ For an ANR spaces $E$ and $E_{2}$, Theorem \ref{Am11} and Corollary \ref{Am12} will remain valid  if we put
   $A$  a closed subspace of metrizable space $X$ instead
of   an ANR pair  $(X,A)$ since we can use Theorem \ref{Am5} in
the proofs. For a polyhedron base $O$,  Corollary \ref{Am12} will
lead us to Theorem \ref{Am7}.}
\end{rem}

\begin{defn}\label{Am14}
\emph{ Let $f:E  \longrightarrow O$ be a   fibration with fiber
space $F_{r_{o}}=f^{-1}(r_{o})$, where $r_{o}\in O$.  By the
\emph{$Lf-$function for  fibration $f$ induced by a lifting
function $L_{f}$ } we mean a map $\Theta_{L_{f}}:
\Omega(O,r_{o})\times F_{r_{o}} \longrightarrow  F_{r_{o}} $ which
is defined by
\[\Theta_{L_{f}}(\alpha,e)=L_{f}(e,\alpha)(1)\quad \mbox{for } \mbox{ }e\in F_{r_{o}}, \alpha\in \Omega(O,r_{o}).\]
}
\end{defn}

Henceforth, we will denote by $[E ,f,O, F_{r_{o}}]$
 the  regular fibration  $f:E  \longrightarrow O$ with   an  $Lf-$function
 $\Theta_{L_{f}}:
  \Omega(O,r_{o})\times F_{r_{o}} \longrightarrow  F_{r_{o}}$, induced by the
 lifting function $L_{f} $ and with a  fiber space
$F_{r_{o}}=f^{-1}(r_{o})$, where  $r_{o}\in O$.


\begin{thm}\label{Am15}
Let $[E ,f,O, F_{r_{o}}]$ be      a fibration with
  metrizable spaces  $E $ and  $O$. Let $\Theta:\Omega(O,r_{o})\times
F_{r_{o}}\longrightarrow F_{r_{o}}$ be
   a map such that
$\Theta_{L_{f}}\simeq \Theta $ and $\Theta(\widetilde{r_{o}},e)=e$
for all $e\in F_{r_{o}}$. If $E$ is an ANR, then there is  a
regular lifting function $L'_{f}$ for $f$ which induces $\Theta$.
That is, $\Theta$ is the $Lf-$function for $f$.
\end{thm}
\noindent \textbf{Proof.}
 Firstly, by the hypothesis, $\Theta_{L_{f}}
\simeq  \Theta$. Then   there is     a homotopy
\[R: [\Omega(O,r_{o})\times F_{r_{o}}]\times I \longrightarrow
 F_{r_{o}} \] such that
$ R[(\alpha,e),0]=\Theta_{L_{f}}(\alpha,e)$ and $
 R[(\alpha,e),1]=\Theta(\alpha,e)$  for all $e\in F_{r_{o}},
 \alpha\in \Omega(O,r_{o})$. We observe that  $\Theta_{L_{f}}$ is extendable to a
  map
$\Theta'_{L_{f}}:  Pa(O,r_{o})\times F_{r_{o}} \longrightarrow E$
defined by
\[\Theta'_{L_{f}}(\alpha,e)=L_{f}(e,\alpha)(1) \quad\mbox{for  }
   \mbox{ } e\in F_{r_{o}},\alpha\in Pa(O,r_{o}) \] having the
   property
\[f[\Theta'_{L_{f}}(\alpha,e)]=\alpha(1) \quad\mbox{for  }
   \mbox{ } e\in F_{r_{o}},\alpha\in Pa(O,r_{o}).\]
This implies that $R$ and $\Theta'_{L_{f}}$ give us a map from
\[[Pa(O,r_{o})\times F_{r_{o}}\times \{0\}]\cup [\Omega(O,r_{o})\times F_{r_{o}}\times I]\] in to $E$.
Since $ \Omega(O,r_{o})\times F_{r_{o}}$ is a closed in a
metrizable space $Pa(O,r_{o})\times F_{r_{o}}$ then by Corollary
\ref{Am12} and Remark \ref{Am13}, $\Theta$ can be extended to a
map
   a map
$\Theta': Pa(O,r_{o})\times F_{r_{o}} \longrightarrow E$ such that
\[f[\Theta'(\alpha,e)]=\alpha(1) \quad\mbox{for  }
   \mbox{ } e\in F_{r_{o}},\alpha\in \Omega(O,r_{o}).\]

Secondly, for $\alpha\in O^I$ and $r\in I$,  define two paths
$\alpha_{r},\alpha^{r}\in O^I$ by
\[\alpha_{r}(t)=\alpha(rt)\quad \mbox{and}\quad \alpha^{r}(t)=\alpha(r+(1-r)t)  \quad\mbox{for   }
   \mbox{ }t\in I.\]
 Hence  define     a homotopy
$H':[\Omega(O,r_{o})\times F_{r_{o}}]\times I\longrightarrow
F_{r_{o}} $ by
\[H'[(\alpha,e),t]=L_{f}[\Theta'(\alpha_{t},e),\alpha^{t}](1) \quad \mbox{for   } \mbox{ }  t\in I, e\in
F_{r_{o}},
 \alpha \in \Omega(O,r_{o}).\]
 Hence
by the hypothesis  and   the  regularity for $L_{f}$, we observe
that
   \begin{eqnarray*}
      H'[(\alpha,e),0]  &=&L_{f}[\Theta'(\alpha_{0},e),\alpha^{0}](1) \\
       &=& L_{f}[\Theta'(\widetilde{r_{o}},e),\alpha](1) \\
       &=& L_{f}(e,\alpha)(1)=\Theta_{L_{f}}(\alpha,e),\\
      H'[(\alpha,e),1] &=&L_{f}[\Theta'(\alpha_{1},e),\alpha^{1}](1) \\
       &=& L_{f}[\Theta'(\alpha,e),\widetilde{r_{o}}](1) \\
       &=&\Theta'(\alpha,e) = \Theta(\alpha,e),
    \end{eqnarray*}
 for all $ e\in
F_{r_{o}}, \alpha \in \Omega(O,r_{o})$ and
\begin{eqnarray*}
      H'[(\widetilde{r_{o}},e),t] &=&L_{f}[\Theta'((\widetilde{r_{o}})_{t},e),(\widetilde{r_{o}})^{t}](1) \qquad \qquad\qquad\quad\mbox{ }\\
       &=& L_{f}[\Theta'(\widetilde{r_{o}},e),\widetilde{r_{o}}](1) \\
       &=& L_{f}(e,\widetilde{r_{o}})(1)= e,
    \end{eqnarray*}
    for all $e\in
F_{r_{o}}$.  Again we can apply  Theorem \ref{Am11} and Remark
\ref{Am13} by taking \[A:=B\cup C\quad \mbox{and}\quad X:=
\bigtriangleup'f=\{(\alpha,e)\in Pa(O)\times E :
 \alpha(0)=f(e)\}\] where $B=[\Omega(O,r_{o})\times
F_{r_{o}}]$, and $C= [(\widetilde{O}\times
 E )\cap \bigtriangleup'f]$,
 and a map $G:(X\times \{0\})\cup(A\times I)\longrightarrow E$
 given by
\[G[(\alpha,e),t]=\left\{%
\begin{array}{ll}
    H'[(\alpha,e),t]  & \mbox{for} \quad [(\alpha,e),t]\in B\times I, \\
    e  & \mbox{for} \quad [(\alpha,e),t]\in C\times I, \\
    L_{f}(e,\alpha)(1)  & \mbox{for} \quad [(\alpha,e),t]\in \bigtriangleup'f\times \{0\}, \\
\end{array}%
\right.\]
      Hence   there is     a map $H: \bigtriangleup'f \times I\longrightarrow E $
     such that $H$ is an extension of $G$ and
     \[f\{H[(\alpha,e),t]\}=f\{H[(\alpha,e),0] \} \quad\mbox{for   }
   \mbox{ } (\alpha,e)\in \bigtriangleup'f, t\in I.\]

   Finally, we can define     a map $L'_{f}: \bigtriangleup f \longrightarrow E^I$ by
    \[L'_{f}(e,\alpha)(t)=H[(\alpha_{t},e),1] \quad\mbox{for   }
   \mbox{ }(e,\alpha)\in \bigtriangleup f,t\in I.\]
    Now we will show that $L'_{f}$ is      a regular lifting function for $f$ as follows:
    \newline
        1.   For $(e,\alpha)\in \bigtriangleup f$, we have that
       \[          L'_{f}(e,\alpha)(0) = H[(\alpha_{0},e),1]=
          G[(\widetilde{\alpha(0)},e),1]=e;\]
        2.  For $(e,\alpha)\in \bigtriangleup f$ and $t\in I$, we have that
        \begin{eqnarray*}
          f[L'_{f}(e,\alpha)(t)] = f\{H[(\alpha_{t},e),1]\}=f\{H[(\alpha_{t},e),0]\}&=&f[L_{f}(e,\alpha_{t})(1)] =\alpha(t);
        \end{eqnarray*}
       3.  For $e\in S$,
        \[          L'_{f}(e,f\circ\widetilde{e})(t) =
        H[(\widetilde{f(e)}_{t},e),1]
           = G[(\widetilde{f(e)}_{t},e),1]=e.\]
    Hence $L'_{f}$ is      a regular lifting function of  $f$ and  for $(\alpha,e)\in \Omega(O,r_{o})\times F_{r_{o}}$,
     \begin{eqnarray*}
       L'_{f}(e,\alpha)(1)= H[(\alpha_{1},e),1]= H[(\alpha,e),1]= H'[(\alpha,e),1]= \Theta(\alpha,e).
     \end{eqnarray*}
    That is,  $\Theta$ is     an $Lf-$function for $f$ induced by the  regular lifting function
     $L'_{f}$.\quad $\square$

It is clear that the lifting function for any fibration no need to
be unique and the definition of the  $Lf-$function depends on the
lifting function. So the $Lf-$function no need to be unique but
 it is uniquely determined up to a
 homotopy class as it is shown in the following theorem.

\begin{thm}\label{Am16}
  Let $[E ,f,O, F_{r_{o}}]$ be     a
 fibration. If $f$ has two  lifting functions $L_{f}$ and
$L'_{f}$, then the $Lf-$functions $\Theta_{L_{f}}$ and
$\Theta_{L'_{f}}$ are homotopic.
\end{thm}
\noindent \textbf{Proof.}
 For $\alpha\in O^I$ and $r\in I$, we can
define two paths $\alpha_{r}$ and $\alpha^{r}$ in $O^I$ by
\[\alpha_{r}(t)=\alpha(rt)\quad\mbox{and}\quad\alpha^{r}(t)=\alpha(r+(1-r)t)\quad\mbox{for}   \mbox{ }t\in I.\] Hence we can define a homotopy
\[H:[\Omega(O,r_{o})\times F_{r_{o}}]\times I\longrightarrow F_{r_{o}} \] by
\[H[(\alpha,e),t]=L_{f}[L'_{f}(e,\alpha_{t})(1),\alpha^{t}](1)\quad\mbox{for  }\mbox{ } t\in I, e\in F_{r_{o}},\alpha\in \Omega(O,r_{o}).\] By
the  regularity for $L_{f}$ and $L'_{f}$, we get that
 \begin{eqnarray*}
     H[(\alpha,e),0]=L_{f}[L'_{f}(e,\alpha_{0})(1),\alpha^{0}](1) &=&
    L_{f}[L'_{f}(e,\widetilde{r_{o}})(1),\alpha](1) \\
       &=& L_{f}(e,\alpha)(1)\\
       &=& \Theta_{L_{f}}(\alpha,e),
       \end{eqnarray*}
       and
\begin{eqnarray*}
     H[(\alpha,e),1]=L_{f}[L'_{f}(e,\alpha_{1})(1),\alpha^{1}](1) &=&
    L_{f}[L'_{f}(e,\alpha)(1),\widetilde{r_{o}}](1)\\
       &=& L'_{f}(e,\alpha)(1)\\
       &=& \Theta_{L'_{f}}(\alpha,e),
\end{eqnarray*}
 for all $ e\in F_{r_{o}}$, $\alpha\in \Omega(O,r_{o})$. Hence $\Theta_{L_{f}}$ and $\Theta_{L'_{f}}$ are
homotopic.\quad $\square$


\begin{defn}\label{Am17}
\emph{Let  $[E_1,f_1,O,F^1_{r_{o}}]$ and $[E_2,f_2,O,F^2_{r_{o}}]$
be two  fibrations. The $Lf-$functions $\Theta_{L_{f_{1}}}$ and
$\Theta_{L_{f_{2}}}$ are said to be \emph{conjugate }if there is
$g\in
 H(F^{1}_{r_{o}},F^{2}_{r_{o}})$ such that
\[\Theta_{L_{f_{1}}}\simeq \overleftarrow{g}\circ
\Theta_{L_{f_{2}}}\circ(i_{\Omega(O,r_{o})}\times g),\] where
$i_{\Omega(O,r_{o})}$ dentes the identity map of $\Omega(O,r_{o})$
onto itself.  }
\end{defn}


\begin{thm}\label{Am18}
Let $[E_1,f_1,O,F^1_{r_{o}}]$ and $[E_2,f_2,O,F^2_{r_{o}}]$ be two
fibrations. If $f_{1}$ and $f_{2}$ are fiber homotopic equivalent,
then $\Theta_{L_{f_{1}}}$ and $\Theta_{L_{f_{2}}}$ are conjugate
$Lf-$functions.
\end{thm}
\noindent \textbf{Proof.}
 Let $f_{1}$ and $f_{2}$ are fiber homotopic equivalent  by two fiber maps  \[h:E_{1}\longrightarrow E_{2}\quad\mbox{and}\quad g:E_{2}\longrightarrow
 E_{1}.\]
For $\alpha\in O^I$ and $r\in I$, we can define two paths
$\alpha_{r}$ and $\alpha'_{r}$ by
\[\alpha_{r}(t)=\alpha(rt)\quad\mbox{and}\quad\alpha'_{r}(t)=\alpha[r+(1-r)t] \quad\mbox{ for }
   \mbox{ } t\in I.\]
Hence we can define      a homotopy $H: \bigtriangleup f_{1}\times
I\longrightarrow  E_{1}$ by
\[H[(\alpha,e),t]=g\bigl\{L_{f_{2}}\{h[L_{f_{1}}(e,\alpha_{t})(1)],\alpha'_{t}\}(1)\bigr\} \quad\mbox{ for }
   \mbox{ }t\in I, (e,\alpha)\in \bigtriangleup f_{1}\] By the
  regularity for $L_{f_{1}}$ and $L_{f_{2}}$, we get that
 \begin{eqnarray*}
     H[(\alpha,e),0]&=&g\bigl\{L_{f_{2}}\{h[L_{f_{1}}(e,\alpha_{0})(1)],\alpha'_{0}\}(1)\bigr\}\qquad\qquad\qquad\quad\\
      &=&g\bigl\{L_{f_{2}}\{h[L_{f_{1}}(e,f_{1}\circ \widetilde{e})(1)],\alpha\}(1)\bigr\} \\
       &=&g[L_{f_{2}}(h(e),\alpha)(1)],
 \end{eqnarray*}
 and
 \begin{eqnarray*}
      H[(\alpha,e),1]&=&g\bigl\{L_{f_{2}}\{h[L_{f_{1}}(e,\alpha_{1})(1)],\alpha'_{1}\}(1)\bigr\}\\
      &=&    g\bigl\{L_{f_{2}}\{h[L_{f_{1}}(e,\alpha)(1)],\widetilde{\alpha(1)}\}(1)\bigr\} \\
       &=&g\bigl\{L_{f_{2}}\{h[L_{f_{1}}(e,\alpha)(1)],f_{1}[\widetilde{L_{f_{1}}(e,\alpha)(1)}]\}(1)\bigr\} \\
      &=& g\bigl\{L_{f_{2}}\{h[L_{f_{1}}(e,\alpha)(1)],(f_{2}\circ h)[\widetilde{L_{f_{1}}(e,\alpha)(1)}]\}(1)\bigr\} \\
      &=& g\bigl\{L_{f_{2}}\{h[L_{f_{1}}(e,\alpha)(1)],f_{2} [\widetilde{h(L_{f_{1}}(e,\alpha)(1))}]\}(1)\bigr\} \\
      &=&g\{h[L_{f_{1}}(e,\alpha)(1)]\}\\
      &=& (g\circ h)[L_{f_{1}}(e,\alpha)(1)],
\end{eqnarray*}
for all $ (e,\alpha )\in \bigtriangleup f_{1}$.
 Consider
the composition
\begin{displaymath} \xymatrix{
 \Omega(O,r_{o})\times F^{1}_{r_{o}}
 \ar@{->}[rr]^{i_{\Omega(O,r_{o})}\times h_{o}} &&  \Omega(O,r_{o})\times
F^{2}_{r_{o}} \ar@{->}[rr]^{\qquad\Theta_{L_{f_{2}}}}&&
F^{2}_{r_{o}} \ar@{->}[r]^{g_{o}}& F^{1}_{r_{o}}}
\end{displaymath}
where $h_{o}=h|_{ F^{1}_{r_{o}} }$ and $g_{o}=g|_{ F^{2}_{r_{o}}
}$. Hence define a homotopy $G:\Omega(O,r_{o})\times
F^{1}_{r_{o}}\times I\longrightarrow  F^{1}_{r_{o}}$ as a
restriction  map  $G=H|_{\Omega(O,r_{o})\times F^{1}_{r_{o}}}$ of
$H$ on $\Omega(O,r_{o})\times F^{1}_{r_{o}}$. Then we get  that
\begin{eqnarray*}
     G[(\alpha,e),0]&=&g_{o}[L_{f_{2}}(h_{o}(e),\alpha)(1)]\\
      &=&    g_{o}[\Theta_{L_{f_{2}}}(\alpha,h_{o}(e)] \\
       &=&[g_{o}\circ \Theta_{L_{f_{2}}}\circ (i_{\Omega(O,r_{o})}\times
       h_{o})](\alpha,e),
       \end{eqnarray*}
       and
       \begin{eqnarray*}
           G[(\alpha,e),1]&=&(g_{o}\circ h_{o})[L_{f_{1}}(e,\alpha)(1)]\qquad\qquad\quad\\
      &=& (g_{o}\circ h_{o})[\Theta_{L_{f_{1}}}(\alpha,e)] \\
       &=&[(g_{o}\circ h_{o})\circ \Theta_{L_{f_{1}}}](\alpha,e),
        \end{eqnarray*}
        for all $e\in F^{1}_{r_{o}}$ and $\alpha\in \Omega(O,r_{o})$.
        Hence
        \[g_{o}\circ \Theta_{L_{f_{2}}}\circ (id_{\Omega(O,r_{o})}\times
        h_{o})\simeq (g_{o}\circ h_{o})\circ
        \Theta_{L_{f_{1}}}.\]
        Since $g_{o}\circ h_{o}\simeq  id_{F^{1}_{r_{o}}}$,
        then $\overleftarrow{h_{o}}\circ \Theta_{L_{f_{2}}}\circ (id_{\Omega(O,r_{o})}\times
        h_{o})\simeq
        \Theta_{L_{f_{1}}}$. Hence $\Theta_{L_{f_{1}}}$
and $\Theta_{L_{f_{2}}}$ are conjugate $Lf-$functions.\quad
$\square$


\begin{lem}\label{Am19}
Let $[E ,f,O, F_{r_{o}}]$ be      a fibration. Then the maps
$D,D_{o}:\bigtriangleup f \longrightarrow E$ defined by
\[D(e,\alpha)=L_{f}[L_{f}(e,\alpha)(1),\overline{\alpha}](1)\quad\mbox{and}\quad D_{o}(e,\alpha)=e,\]
for all $(e,\alpha)\in \bigtriangleup f$,
 are  homotopic.
 \end{lem}
\noindent \textbf{Proof.}
 For $\alpha \in O^I$ and $r\in I$,   define  paths
$\alpha_{r},\alpha'_{r}$ and $ \alpha''_{r}$ in $O$ by
\[\alpha_{r}(t)=\alpha(rt),\quad\alpha'_{r}(t)=\alpha[r+(1-r)t] \quad\mbox{and}\quad \alpha''_{r}(t)=\alpha[2r(1-t)],\]
 for all $t\in I$.
Define     two  homotopies $H: \bigtriangleup f \times
I\longrightarrow E$ by
\[H[(e,\alpha),t]=L_{f}[L_{f}(e,\alpha_{t})(1),\alpha'_{t}](1) \quad\mbox{ for}
   \mbox{ }t\in I, (e,\alpha)\in \bigtriangleup f,\] and a homotopy
   $G:O^I \times I\longrightarrow O^I$ by
\[ [G(\alpha,r)](t) = \left\{
  \begin{array}{c l}
    \alpha_{r}(t) & \mbox{for \quad$ 0 \leq t \leq 1/2$},\\
    \alpha''_{r}(t)  & \mbox{for\quad $1/2 \leq t \leq 1$},
  \end{array}
\right. \] for all  $\alpha\in O^I, r\in I$.
   Hence define a homotopy $F: \bigtriangleup f \times I\longrightarrow
E$ by
\[F[(e,\alpha),t]=H[(e,G(\alpha,t)),1/2] \quad\mbox{ for}
   \mbox{ } t\in I, (e,\alpha)\in \bigtriangleup f.\]
   By the  regularity for $L_{f}$     we observe that
for $(e,\alpha)\in \bigtriangleup f$,
\begin{eqnarray*}
F[(e,\alpha),1]=H[(e,G(\alpha,1)),1/2]&=&H[(e,\alpha\star\overline{\alpha}),1/2]\\
   &=& L_{f}\{K[e,(\alpha\star\overline{\alpha})_{1/2}],(\alpha\star\overline{\alpha})'_{1/2}\}(1) \\
   &=& L_{f}[K(e,\alpha),\overline{\alpha}](1) \\
   &=& L_{f}[L_{f}(e,\alpha)(1),\overline{\alpha}](1)  \\
   &=& D(e,\alpha)
\end{eqnarray*}
  for all $(e,\alpha)\in \bigtriangleup f$, and
  \begin{eqnarray*}
F[(e,\alpha),0]=H[(e,G(\alpha)),0](1/2)&=&H[(e,\widetilde{\alpha(0)}),1/2]\\
   &=& L_{f}\{K[e,\widetilde{\alpha(0)}_{1/2}],\widetilde{\alpha(0)}'_{1/2}\}(1) \\
   &=& L_{f}\{K[e,\widetilde{\alpha(0)}],\widetilde{\alpha(0)}\}(1) \\
   &=& L_{f}[L_{f}(e,f\circ\widetilde{e})(1),f\circ\widetilde{e}](1)  \\
   &=& L_{f}(e,f\circ\widetilde{e})(1)\\
   &=& e= D_{o}(e,\alpha)
\end{eqnarray*} for all $(e,\alpha)\in \bigtriangleup f$.
Hence $D$ and $D_{o}$ are  homotopic.\quad $\square$

 In the proof
of Lemma  above we get that the  homotopy $F$ has the following
property:\newline
\begin{eqnarray}\label{eeq3}
f\{F[(e,\alpha),t]\}   &=&\alpha(0) \quad\mbox{ for}
   \mbox{ } (e,\alpha)\in \bigtriangleup f.
   \end{eqnarray}

\begin{rem}\label{Am20}
\emph{In  Lemma  \ref{Am19},  for $t\in I$, we can define a map
$\Theta^{t}_{L_{f}}: \Omega(O,r_{o})\times F_{r_{o}}
\longrightarrow F_{r_{o}} $ by
\[\Theta^{t}_{L_{f}}(\alpha,s)=H(s,\alpha)(t) \quad\mbox{ for  }
   \mbox{ }\alpha\in
\Omega(O,r_{o}),s\in F_{r_{o}}.\] By the   regularity for $L_{f}$,
we observe that
\[  \Theta^{0}_{L_{f}}(\alpha,s) = \Theta^{1}_{L_{f}}(\alpha,s) =
\Theta_{L_{f}}(\alpha,s).\]
  By the
definition of $\Theta^{t}_{L_{f}}$ we get that for $t\in I$,
$\Theta^{t}_{L_{f}}\simeq \Theta_{L_{f}}$. Also we observe that
$\Theta^{t}_{L_{f}}(\widetilde{r_{o}},s)=s$ for all $s\in
F_{r_{o}}$.}
\end{rem}

\noindent
 \noindent \textbf{Proof of Theorem (1.2).}
 \emph{Necessity}: If $f_{1}$ and $f_{2}$ are fiber homotopy
equivalent then by Theorem \ref{Am18} they have conjugate
$Lf-$functions.

 \emph{Sufficiency}: Since $O_{1}$ and $O_{2}$
are   contractible to   $r_{o}\in O_{3}$, then there are two
 homotopy maps $R_{1}:O_{1}\times I\longrightarrow O$ and
$R_{2}:O_{2}\times I\longrightarrow O$ such that
\[R_{1}(x,0)=x \quad R_{1}(x,1)=r_{o} \quad\mbox{ for }
   \mbox{ }x\in
O_{1},\]and
\[R_{2}(x,0)=x \quad R_{2}(x,1)=r_{o} \quad\mbox{ for }
   \mbox{ }x\in
O_{2},\]respectively. We will denote  the path: $t\longrightarrow
R_{1}(x,t)$ by $R^{x}_{1}$ for $x\in O_{1}$ and the path:
$t\longrightarrow R_{2}(x,t)$ by $R^{x}_{2}$ for $x\in O_{2}$.

 In Figure 1, let  $f_{ij}=f_{i}|O_{j}$, where $i=1,2$ and $j=1,2,3 $. Now  we can define
  a map $h_{1}: f^{-1}_{1}(O_{1}) \longrightarrow
 f^{-1}_{2}(O_{1}) $  by
\[h_{1}(e)=L_{f_{2}}\{g[L_{f_{1}}(e,R^{f_{1}(e)}_{1})(1)],\overline{R^{f_{1}(e)}_{1}}\}(1) \quad\mbox{ for  }
   \mbox{ }e\in f^{-1}_{1}(O_{1}), \] and      a map
$h_{2}: f^{-1}_{1}(O_{2}) \longrightarrow f^{-1}_{2}(O_{2}) $ by
\[h_{2}(e)=L_{f_{2}}\{g[L_{f_{1}}(e,R^{f_{1}(e)}_{2})(1)],\overline{R^{f_{1}(e)}_{2}}\}(1) \quad\mbox{ for  }
   \mbox{ }e\in f^{-1}_{1}(O_{2}).\] It is   clear that
\[(f_{2}\circ
h_{1})(e)=\overline{R^{f_{1}(e)}_{1}}(1)=R^{f_{1}(e)}_{1}(0)=f_{1}(e)
\] for all $e\in  f^{-1}_{1}(O_{1})$ and
\[(f_{2}\circ
h_{2})(e)=\overline{R^{f_{1}(e)}_{2}}(1)=R^{f_{1}(e)}_{2}(0)=f_{1}(e)
\] for all $e\in  f^{-1}_{1}(O_{2})$. That is, $h_{1}$ and $h_{2}$ are  fiber maps.
   \begin{figure}[h!]
   \begin{center}
  \[ \xymatrix{
O_1                                && f^{-1}_1(O_1)  \ar[ll]_{f_{11}} \ar[rr]^{h_1}                                                                                        && f^{-1}_2(O_1) \ar[rr]^{f_{21}}                            && O_1 \\
 O_3  \ar@{.>}[u]^{i} \ar@{.>}[d]_{i}  && f^{-1}_1(O_3)  \ar@{.>}[u]^{i}\ar@{.>}[d]_{i}\ar[ll]_{f_{13}} \ar@{->}   @<4pt>[rr]^{h_1\circ k_1}\ar@{->}   @<-4pt>[rr]_{k_2 \circ h_2} && f^{-1}_2(O_3) \ar@{.>}[u]^{i}\ar@{.>}[d]_{i}\ar[rr]^{f_{23}}  && O_3  \ar@{.>}[u]^{i}\ar@{.>}[d]_{i} \\
 O_2                               && f^{-1}_1(O_2)  \ar[ll]_{f_{12}} \ar[rr]^{h_2}                                                                                        && f^{-1}_2(O_2) \ar[rr]^{f_{22}}                            && O_2
 }
\]
  \end{center}
  \vspace{-5mm}
 \center{Figure 1}
\end{figure}

 Also in Figure 1, we can define      a map
$k_{1}:f^{-1}_{1}(O_{3}) \longrightarrow f^{-1}_{1}(O_{3}) $ by
\[k_{1}(e)=L_{f_{1}}[L_{f_{1}}(e,R^{f_{1}(e)}_{2})(1),\overline{R^{f_{1}(e)}_{2}}](1) \quad\mbox{ for   }
   \mbox{ }e\in f^{-1}_{1}(O_{3}), \] and      a map
$k_{2}:f^{-1}_{2}(O_{3}) \longrightarrow f^{-1}_{2}(O_{3}) $ by
\[k_{2}(e)=L_{f_{2}}[L_{f_{2}}(e,R^{f_{2}(e)}_{1})(1),\overline{R^{f_{2}(e)}_{1}}](1) \quad\mbox{ for   }
   \mbox{ }e\in f^{-1}_{2}(O_{3}),\] respectively. It is   clear that
\[(f_{1}\circ
k_{1})(e)=\overline{R^{f_{1}(e)}_{2}}(1)=R^{f_{1}(e)}_{2}(0)=f_{1}(e)
 \] for all $e\in f^{-1}_{1}(O_{3})$ and
\[(f_{2}\circ
k_{2})(e)=\overline{R^{f_{2}(e)}_{1}}(1)=R^{f_{2}(e)}_{1}(0)=f_{2}(e)
\] for all $e\in f^{-1}_{2}(O_{3})$.  That is, $k_{1}$ and $k_{2}$ are fiber maps. By Lemma
\ref{Am19}, we get that
\[k_{1}\simeq_{f}i_{f^{-1}_{1}(O_{3})}\quad\mbox{and}\quad
k_{2}\simeq_{f}i_{f^{-1}_{2}(O_{3})}.\]Then
\[h_{1}\circ
k_{1}\simeq_{f}h_{1} \quad\mbox{and}\quad k_{2}\circ
h_{2}\simeq_{f}h_{2}.\] That is, to show that
$h_{1}\simeq_{f}h_{2}$ as  fiber maps of $f^{-1}_{1}(O_{3}) $ into
$f^{-1}_{2}(O_{3}) $, it is sufficient to show that $h_{1}\circ
k_{1}\frac{f}{\simeq}k_{2}\circ h_{2}$ as  fiber maps of
$f^{-1}_{1}(O_{3}) $ into $f^{-1}_{2}(O_{3}) $. By Remark
\ref{Am20}, we get that
\begin{eqnarray*}
  (h_{1}\circ k_{1})(e) &=& L_{f_{2}}\bigl\{g\{\Theta^{1/2}_{L_{f_{1}}}[\overline{R^{f_{1}(e)}_{2}}\star
R^{f_{1}(e)}_{1},L_{f_{1}}(e,R^{f_{1}(e)}_{2})(1)]\},\overline{R^{f_{1}(e)}_{1}}\bigr\}(1) \\
   &=& L_{f_{2}}\{(g\circ \Theta^{1/2}_{L_{f_{1}}})[\overline{R^{f_{1}(e)}_{2}}\star
R^{f_{1}(e)}_{1},L_{f_{1}}(e,R^{f_{1}(e)}_{2})(1)],\overline{R^{f_{1}(e)}_{1}}\}(1),
\end{eqnarray*}
and
\begin{eqnarray*}
  (k_{2}\circ h_{2})(e)&=& L_{f_{2}}\bigl\{\Theta^{1/2}_{L_{f_{2}}}\{\overline{R^{f_{1}(e)}_{2}}\star
R^{f_{1}(e)}_{1},g[L_{f_{1}}(e,R^{f_{1}(e)}_{2})(1)]\},\overline{R^{f_{1}(e)}_{1}}\bigr\}(1) \\
   &=& L_{f_{2}}\bigl\{[\Theta^{1/2}_{L_{f_{2}}}\circ(1\times g)]\{\overline{R^{f_{1}(e)}_{2}}\star
R^{f_{1}(e)}_{1},L_{f_{1}}(e,R^{f_{1}(e)}_{2})(1)\},\overline{R^{f_{1}(e)}_{1}}\bigr\}(1),
\end{eqnarray*}
for all $e\in f^{-1}_{1}(O_{3})$, where $1=i_{\Omega(O,r_{o})}$.

By the hypothesis $\Theta_{L_{f_{1}}}\simeq \overleftarrow{g}\circ
\Theta_{L_{f_{2}}}\circ(i_{\Omega(O,r_{o})}\times g)$, e.g.,
\[g\circ \Theta_{L_{f_{1}}}\simeq
\Theta_{L_{f_{2}}}\circ(i_{\Omega(O,r_{o})}\times g),\] and by
Remark \ref{Am20}, we get $\Theta_{L_{f_{1}}}\simeq
\Theta^{1/2}_{L_{f_{1}}}$ and $ \Theta_{L_{f_{2}}}\simeq
\Theta^{1/2}_{L_{f_{2}}}$. That is,
\begin{eqnarray*}
  g\circ \Theta^{1/2}_{L_{f_{1}}} &\simeq & \Theta^{1/2}_{L_{f_{2}}}\circ(i_{\Omega(O,r_{o})}\times
  g).
\end{eqnarray*} Hence  by
 Lemma \ref{Am19} again, we get
\begin{eqnarray}\label{eq}
  h_{1}\circ
k_{1}\simeq_{f}k_{2}\circ h_{2}\quad &\Longrightarrow &\quad
h_{1}\simeq_{f}h_{2},
\end{eqnarray}
 as  fiber maps of
$f^{-1}_{1}(O_{3}) $ into $f^{-1}_{2}(O_{3}) $.

Now in two fibrations   $f_{13}=f_{1}|O_{3}$ and
$f_{23}=f_{2}|O_{3}$, since $O_{3}$ is a subpolyhedra of $O_{2}$,
$h_{1}\simeq_{f}h_{2}$ as maps of $f^{-1}_{1}(O_{3}) $ into
$f^{-1}_{2}(O_{3}) $ and $h_{2}$ is defined on
$f^{-1}_{1}(O_{2})$, then by Theorem \ref{Am7} $h_{1}$ can be
extended as      fiber map to all of $f^{-1}_{1}(O_{2})$. Since
$h_{1}$ is  defined on $ f^{-1}_{1}(O_{1})$ then $h_{1}$ gives a
fiber map $h$ of $ E_{1}$ into $ E_{2} $.
 Since $O_{1}$ is     contractible to $r_{o}$ leaves $r_{o}$ a fixed,  then
$R^{r_{o}}_{1}=\widetilde{r_{o}}$. Hence for $e\in
F^{1}_{r_{o}}\subseteq f^{-1}_{1}(O_{3})$,
\begin{eqnarray*}
  h(e) = h_{1}(e)&=& L_{f_{2}}\{g[L_{f_{1}}(e,R^{f_{1}(e)}_{1})(1)],\overline{R^{f_{1}(e)}_{1}}\}(1)\\
   &=& L_{f_{2}}\{g[L_{f_{1}}(e,R^{r_{o}}_{1})(1)],\overline{R^{r_{o}}_{1}}\}(1)\\
   &=& L_{f_{2}}\{g[L_{f_{1}}(e,\widetilde{r_{o}})(1)],\overline{\widetilde{r_{o}}}\}(1)\\
   &=& L_{f_{2}}(g(e),\widetilde{r_{o}})(1) = g(e).
\end{eqnarray*}
That is, $h$  as     map: $F^{1}_{r_{o}}\longrightarrow
F^{2}_{r_{o}}$ is      a  homotopy equivalent. Since  $O$ is an
ANR and it is clear that $ O $ is a pathwise
 connected ($ O $ is the
union for two contractible spaces and $O_{3}\neq\phi$), then by
Fadell-Dold theorem, $f_{1}$ and $f_{2}$ are fiber homotopy
equivalent.  \quad $\square$


\section{Applying $Lf-$function in fiber bundles}
Here we apply the $Lf-$function in fiber bundles by proving  the
equivalently between Theorem \ref{Am2} and Dold's theorem.

Firstly, we will give some propositions which help us to make
comparing between Dold's theorem and Theorem \ref{Am2}.

In the following proposition we will prove the converse of Theorem
\ref{Am10}.
\begin{pro}\label{Am21}
\emph{Let $G$ be a group of all homeomorphisms of space $F$ with
as binary usual composition operation $\circ$. If there is a map
$\mu:(S^{n-1}, x_o)\longrightarrow (G,\mathbf{g})$, then there is
bundle $E$ over sphere $S^n$ and a map $f:E\longrightarrow S^n$
such that $\gamma=(E,f,S^n,F,G)$ is fiber bundle.}
\end{pro}
\noindent\textbf{Proof.}
 Let $S^n=V_1
\cup V_2$, where each $V_i$ is an open n-cell such that $V_1\cap
V_2$ is a trip a round $S^{n-1}$ and there is a retraction
$r:V_1\cap V_2\longrightarrow S^{n-1}$ ( see\cite{James}).
 Now define maps
 \[g_{ii}:V_i\longrightarrow G, \quad g_{ii}(x)=\mathbf{g}\quad \forall
 \mbox{ }x\in V_i, (i=1,2),\]
 \[g_{12}:V_1\cap V_2\longrightarrow G, \quad g_{12}(x)=(\mu \circ r)(x)\quad \forall
 \mbox{ }x\in V_1\cap V_2,\]and
 \[g_{21}:V_1\cap V_2\longrightarrow G, \quad g_{21}(x)=[g_{12}(x)]^{-1}\quad \forall
 \mbox{ }x\in V_1\cap V_2.\]
 Let $J=\{1,2\}$ be a space with the discrete topology and
 $T\subset S^n \times F\times J$ be the set defined by
 \[T=\{(x,y,j):x\in V_j, y\in F, j\in J\}.\]
 Define an equivalent relation $\equiv$ on $T$ by
 \[(x_1,y_1,j)\equiv (x_2,y_2,k) \Longleftrightarrow x_1=x_2 \quad
 \mbox{and} \quad g_{kj}(x_1)(y_1)=y_2,\] where $(x_1,y_1,j),(x_2,y_2,k)
 \in T$. Then put $E$ to be the space of equivalence classes
 obtained with the quotient topology.
 Hence define a map $f:E\longrightarrow S^n$ by
 \[ f([(x,y,j)])=x \quad \forall
 \mbox{ }[(x,y,j)]\in E,\] and the maps $
\theta_j:V_j \times F\longrightarrow f^{-1}(V_j)$ defined by
\[\theta_j(x,y)=[(x,y,j)]\quad \forall
 \mbox{ } (x,y)\in V_j \times F.\] Hence it is clear that $\gamma=(E,f,S^n,F,G)$ is a fiber
 bundle. \quad$\square$


\begin{pro}\label{Am22}
\emph{Let $[E ,f,O, F_{r_{o}}]$ be a   fibration with locally
compact fiber space $F_{r_{o}}$. Then the function
$\phi:\Omega(O,r_{o})\longrightarrow F_{r_{o}}^{F_{r_{o}}}$ given
by
\[ \phi(w)(e)=\Theta_{L_{f}}(w,e)\quad \mbox{for} \mbox{ } w\in
\Omega(O,r_{o}), e\in F_{r_{o}},\] is a map from $\Omega(O,r_{o})$
into $H(F_{r_{o}},F_{r_{o}})$.}
\end{pro}
\noindent \textbf{Proof.}
 Since a Hausdorff space  $F_{r_{o}}$ is a locally compact  then $\phi$ is continuous function (see \cite{Allen} Proposition A.14  P.530).
 Now we will prove that
for $w\in \Omega(O,r_{o})$, $\phi(w)$ is homotopy equivalence from
$F_{r_{o}}$ into $F_{r_{o}}$. For $w\in \Omega(O,r_{o})$, we can
define a map $\overleftarrow{\phi(w)}:F_{r_{o}}\longrightarrow
F_{r_{o}}$ by
\[
\overleftarrow{\phi(w)}(e)=\Theta_{L_{f}}(\overline{w},e)\quad
\mbox{for} \mbox{ }   e\in F_{r_{o}}.\] Then we get that
\[[\phi(w)\circ
\overleftarrow{\phi(w)}](e)=L_{f}[L_{f}(e,\overline{w})(1),w](1)\quad
\mbox{for} \mbox{ } e\in F_{r_{o}},\] and
\[\overleftarrow{\phi(w)}\circ \phi(w)](e)=L_{f}[L_{f}(e,w)(1),\overline{w}](1)\quad
\mbox{for} \mbox{ } e\in F_{r_{o}}.\] Then by Lemma  \ref{Am19},
\[\phi(w)\circ \overleftarrow{\phi(w)}\simeq i_{F_{r_{o}}} \quad \mbox{and}\quad
\overleftarrow{\phi(w)}\circ \phi(w)\simeq i_{F_{r_{o}}}.\] Hence
 $\phi(w)\in H(F_{r_{o}},F_{r_{o}})$. Therefore
$\phi$ is a map from $\Omega(O,r_{o})$ into
$H(F_{r_{o}},F_{r_{o}})$.\quad $\square$
\begin{pro}\label{Am23}
\emph{Let $\gamma=[E,f,O,F_{r_{o}},G]$ be a fiber bundle and
fibration with locally compact fiber $F_{r_{o}}$. Then the
function $\phi:\Omega(O,r_{o})\longrightarrow
F_{r_{o}}^{F_{r_{o}}}$ given by \[
\phi(w)(e)=\Theta_{L_{f}}(w,e)\quad \mbox{for} \mbox{ } w\in
\Omega(O,r_{o}), e\in F_{r_{o}},\] is a map from $\Omega(O,r_{o})$
into $G$.}
\end{pro}
\noindent \textbf{Proof.}
 Since $F$ is a locally compact, then $\phi$ is continuous function.  For $w\in O^I$, let
$F_{w(0)}=f^{-1}(w(0))$ and   $F_{w(1)}=f^{-1}(w(1))$. Then the
map $h: F_{w(0)}\longrightarrow F_{w(1)}$ given by
\[h(e)=L_{f}(e,w)(1)\quad \mbox{for} \mbox{ }e\in  F_{w(0)},\] is a homeomorphism since it is obtained
from the compositions of coordinate functions which are
homeomorphisms. Hence  $\phi$ is a map from $\Omega(O,r_{o})$ into
$G$.\quad $\square$
\begin{defn}\label{Am24}
\emph{For any space $O$ with fixed point $r_o\in O$,  we can
define  a \emph{conical map}  $\psi:O\longrightarrow
\Omega(S(O),r_o)$ as follows:\newline
          For $x\in O$, let $w_0(x)$ be path between in a cone $S_0(O)$ from  equivalent class $[(x,1/2)]$ into $[(r_o,1/2)]$
       and let $w_1(x)$ be path between in a cone $S_1(O)$ from $[(x,1/2)]$ into $[(r_o,1/2)]$. Define the \emph{conical map }$\psi:O\longrightarrow \Omega(S(O),r_o)$ by
       \[\psi(x)=\overline{w_1(x)}  \star w_0(x)\quad \mbox{for} \mbox{ }x\in  O,\]
       where $\Omega(S(O),r_o) := \Omega(S(O),[(r_o,1/2)])$.}
\end{defn}

To prove the equivalently between Theorem  \ref{Am2} and the
Dold's theorem we will  rephrase  Theorem \ref{Am2} for two
  fibrations over a common suspension base.


\begin{rem}\label{Am25}
\emph{ Let $\gamma_1 =[E_{1} ,f_{1} ,S(O),F^{1}_{r_{o}} ,G_{1} ] $
and $\gamma_2=[E_{2},f_{2},S(O),F^{2}_{r_{o}},G_{2}]$ be two
fibrations over a common suspension base $S(O)$  of a polyhedron
space $O$ with locally compact fibers $F^{1}_{r_{o}} $and
$F^{2}_{r_{o}}$. In Figure 2, let
\[\mu_{1} :(O, r_o)\longrightarrow (G_{1} ,\mathbf{g}_{1} )\quad
\mbox{and}\quad\mu_{2}:(O, r_o)\longrightarrow
(G_{2},\mathbf{g}_{2})\] be characteristic maps of $\gamma $ and
$\gamma'$, respectively. Also let
\[i_{1} :G_{1} \longrightarrow H(F^{1}_{r_{o}} , F^{1}_{r_{o}} )\quad \mbox{and}\quad
i_{2}:G_{2}\longrightarrow H(F^{2}_{r_{o}}, F^{2}_{r_{o}})\] be
the inclusion maps. From Propositions  \ref{Am22} and
 \ref{Am23}, then Theorem \ref{Am2} and the Dold's theorem
can now be compared.
\begin{figure}[h!]
   \begin{center}
  \[ \xymatrix{
                                                                                       && L(S(O),r_{o})\\
 G_{2}  \ar[rru]^{\phi_2} \ar[d]^{i_2}  && O  \ar[u]^{\psi} \ar[ll]_{\mu_2}  \ar[rr]^{\mu_1} && G_{1}  \ar[llu]_{\phi_1} \ar[d]^{i_1}\\
 H(F^{2}_{r_{o}} , F^{2}_{r_{o}} )                              &&                                                                                     && H(F^{1}_{r_{o}} , F^{1}_{r_{o}} )\ar[llll]^{T(f)=g\circ f\circ \overleftarrow{g}\mbox{ }\forall \mbox{ } f\in H(F^{1}_{r_{o}} , F^{1}_{r_{o}} ) }
 }
\]
  \end{center}
  \vspace{-5mm}
 \center{Figure 2}
\end{figure}
 Let $g\in H(F^{1}_{r_{o}},F^{2}_{r_{o}})$ and $\psi$ is the conical map.
 Hence  Theorem  \ref{Am2}  can be restated in terms of $\phi_{1}$,
 $\phi_{2}$, and $\psi$ as follows:\newline
 Two   fibrations   $\gamma_1 =[E_{1} ,f_{1} ,S(O),F^{1}_{r_{o}} ,G_{1}]
 $ and $\gamma_2=[E_{2},f_{2},S(O),F^{2}_{r_{o}},G_{2}]$  are   fiber
homotopy equivalent if and only if there is $g\in
H(F^{1}_{r_{o}},F^{2}_{r_{o}})$ such that two maps
\[m(x)=g\circ i_{1}\circ \phi_{1}[\psi(x)]\circ \overleftarrow{g}\quad \mbox{for} \mbox{
} x\in O,\] and
\[m'(x)= i_{2}\circ \phi_{2}[\psi(x)]\quad \mbox{for} \mbox{
} x\in O,\] from $O$ into $H(F^{2}_{r_{o}},F^{2}_{r_{o}})$ are
homotopic.}
\end{rem}

If it can be shown that $\phi_{1} \circ \psi\simeq \mu_{1}$ and
$\phi_{2} \circ \psi\simeq \mu_{2}$, then Theorem  \ref{Am2} and
the Dold's theorem are equivalent. We will prove it in Theorem
\ref{Am27}.

\begin{lem}\label{Am26}
\emph{Let $\gamma=[E,f,O,F_{r_{o}},G]$ be a
  fibration. Let $  O = O_1\cup O_2$, $ r_{o}\in O_1\cap O_2$,
  \begin{eqnarray*}
   \overline{\bigtriangleup}  f &=& \{(\beta,e)\in \Omega(O,r_{o})\times F_{r_{o}} :
\beta=w_2\star w_1, \mbox{ }w_i\in O^I_i
(i=1,2), \\
   &\quad  & \beta(1/2)=w_2(1)=w_1(0)\in O_1\cap
O_2\}
\end{eqnarray*}
 and  $L_{i}$  be a  lifting functions for   fibration   $f|O_{i}$. If there are  fiber homeomorphisms
$\epsilon_i:O_i\times F_{r_{o}}\longrightarrow f^{-1}(O_i)$, then
the map $\overline{\phi}:\overline{\bigtriangleup}
f\longrightarrow F_{r_{o}}$ given by
\[ \overline{\phi}(\beta,e)=L_{1}[L_{2}(e,w_2)(1),w_1](1)
\quad\mbox{for}\mbox{ } (\beta,e)\in \overline{\bigtriangleup}
f,\] is homotopic to the $Lf-$function $\Theta_{L_{f}}$, where
$i=1,2$.}
\end{lem}
\noindent \textbf{Proof.}
 We can define the   lifting functions  $L_{1}$
and $L_{2}$ for   fibrations $f|O_{1}$ and $f|O_{2}$ by
\[L_{i}(e,w)=\epsilon_i[w(t),(\pi_2\circ \epsilon^{-1}_i)(e)]\quad \mbox{for} \mbox{ }  (e,w)\in \bigtriangleup f|O_{i},\]
where $i=1,2$, respectively. Since $L_{f}$ is   lifting function
for $\gamma$, then it is also a   lifting function for $f|O_{1}$
and $f|O_{2}$. Hence $L_{f} \simeq L_{1}$ on  $\bigtriangleup
f|O_{1}$ and $L_{f} \simeq L_{2}$ on  $\bigtriangleup f|O_{2}$.
Define  a map $\widetilde{\phi}:\overline{\bigtriangleup}
f\longrightarrow F_{r_{o}}$
 by
\[ \widetilde{\phi}(\beta,e)=L_{f} [L_{f} (e,w_2)(1),w_1](1)
\quad\mbox{for}\mbox{ } (\beta,e)\in \overline{\bigtriangleup}
f.\] Then  $ \widetilde{\phi}\simeq \overline{\phi}$ and by the
homotopy $H$ in proof of Lemma \ref{Am19}, $
\widetilde{\phi}\simeq \Theta_{L_{f}}$. Hence $
\Theta_{L_{f}}\simeq \overline{\phi}$.\quad $\square$
\begin{thm}\label{Am27}
\emph{Let $\gamma=[E,f,S(O),F_{r_{o}},G]$ be a fiber bundle over
suspension $S(O)$ of a polyhedron space $O$   with locally compact
fiber $F_{r_{o}}$ and admits a lifting function $ L_{f}$. Also let
$\phi:\Omega(S(O).r_{o})\longrightarrow G$ be a map given by
\[ \phi(\beta)(x)=L_{f}(x,\beta)(1)\quad \mbox{for} \mbox{ } \beta\in
\Omega(S(O),r_{o}), x\in F_{r_{o}}.\] Then $\phi \circ \psi\simeq
\mu$, where $\mu:(O,r_{o})\longrightarrow (G,e)$ is the
characteristic map of $\gamma$ and $\psi$ is the conical map.}
\end{thm}
\noindent \textbf{Proof.}
 In Lemma  \ref{Am26}, put $S(O)=S_0(O)\cup S_1(O)$,
$O_1=S_0(O)$, and $O_2=S_1(O)$. It's clear that $O=S_0(O)\cap
S_1(O)$. Now define maps
 \[g_{ii}:S_i(O)\longrightarrow G\quad\mbox{by} \quad g_{ii}(x)=\mathbf{g}\quad \mbox{for}
 \mbox{ }x\in S_i(O), (i=0,1),\]
 \[g_{01}:O\longrightarrow G\quad\mbox{by} \quad g_{01}(x)= \mu(x)\quad \mbox{for}
 \mbox{ }x\in O,\]and
 \[g_{10}:O\longrightarrow G\quad \mbox{by} \quad g_{10}(x)=[\mu(x)]^{-1}\quad \mbox{for}
 \mbox{ }x\in O.\]
 Let $J=\{0,1\}$ be a space with the discrete topology and let
 $T\subset S(O) \times F_{r_{o}}\times J$ be the set defined by
 \[T=\{(x,e,j):x\in S_j(O), e\in F_{r_{o}}, j\in J\}.\]
 Define an equivalent relation $\equiv$ on $T$ by
 \[(x_1,e_1,j)\equiv (x_2,e_2,k) \Longleftrightarrow x_1=x_2 \quad
 \mbox{and} \quad g_{kj}(x_1)(e_1)=e_2,\] where $(x_1,e_1,j),(x_2,e_2,k)
 \in T$.

Recall Proposition \ref{Am21} that points of $E$ are identified to
the equivalent classes of all triples $(x,e,j)\in T$. Hence for
$j=0,1$, the maps $\epsilon_j:S_j(O)\times
F_{r_{o}}\longrightarrow f^{-1}(S_j(O))$ given by
\[\epsilon_j(x,e)=[(x,e,j)]\quad \mbox{for}
 \mbox{ } (x,e)\in S_j(O) \times F_{r_{o}},\] denotes the equivalence class
 of the triple $(x,e,j)$.

  Put $r_{o}:=[(r_{o},1/2)]$ in $S(O)$ and  $e:=[(r_{o},e,i)]$, where $i=0,1$.
 Then from Lemma \ref{Am26}, we have that for $\beta\in \Omega(S(O),r_{o})$,
 $\beta=w_1 \star w_0$  for some $w_0\in [S_{0}(O)]^I, w_1\in [S_{1}(O)]^I$ and
 \begin{eqnarray*}
   L_{1}(e,w_0)(1) &=& [(w_0(1),e,0)], \\
   L_{2}(e,w_1)(1)&=& [(w_1(1),e,1)]\\
     &=&[(w_1(1),\mu(w_1(1))(e),0)].
 \end{eqnarray*}
 Hence
\begin{eqnarray*}
  \overline{\phi}(\beta,e)&=&L_{1}[L_{2}(e,w_1)(1),w_0](1) \\
   &=& [(r_{o},\mu(w_1(1))(e),0)] \\
    &=& \mu(w_1(1))(e) \\
    &=& \mu(\beta(1/2))(e).
  \end{eqnarray*}
  Let $\overline{\Omega}(S(O),r_{o})$ be the projection of
  $\overline{\bigtriangleup}f$ on $\Omega(S(O),r_{o})$ and
  $\phi':\overline{\Omega}(S(O),r_{o})\longrightarrow G$ be a map given by
\[  \phi'(\beta)(e)=\overline{\phi}(\beta,e)\quad\mbox{for}\mbox{ }  \beta \in \overline{\Omega}(S(O),r_{o}), e\in F_{r_{o}}.\]
Then by Lemma \ref{Am19}, $\phi \simeq \phi'$ and
\[(\phi'\circ \psi)(e)=\phi'[\psi(e)]=\mu[\psi(e)(1/2)]=\mu(e),\]
Thus  $\phi'\circ \psi=\mu$. Hence  $\phi\circ \psi\simeq\mu$.
Therefore Theorem  \ref{Am2} and the Dold's theorem are
equivalent.\quad $\square$


\bibliographystyle{plain}

\end{document}